\documentclass[11pt]{article}

\usepackage[a4paper,left=2cm,right=2cm,top=1.5cm,bottom=2cm]{geometry}
\usepackage[utf8]{inputenc}
\usepackage[T1]{fontenc}
\usepackage[english]{babel}
\usepackage{amsthm}
\usepackage{amsmath}
\usepackage{amssymb}
\usepackage{mathrsfs}
\usepackage{setspace}
\usepackage{float}
\usepackage{multirow}
\usepackage{stmaryrd}
\usepackage{color}
\usepackage{wrapfig}
\usepackage{pifont}
\usepackage{array}
\usepackage[colorlinks=true,linkcolor=ocre,citecolor=ocre]{hyperref}
\usepackage{dsfont}
\usepackage{upgreek}
\usepackage{nicefrac}
\usepackage{tikz}{}
\usepackage{gensymb}
\usepackage{FiraSans}

\usetikzlibrary{arrows}
\usetikzlibrary{shapes}
\usetikzlibrary{positioning}
\usetikzlibrary{decorations.pathmorphing}

\usepackage{graphicx}
\usepackage{caption}
\renewenvironment{proof}{\noindent{\sffamily{\textbf{Proof :}}}}{\begin{flushright}$\square$\end{flushright}}

\newcommand{\IE}{\mathbb{E}}

\newcommand{\IZ}{\mathbb{Z}}
\newcommand{\IQ}{\mathbb{Q}}
\newcommand{\IR}{\mathbb{R}}

\newcommand{\IT}{\mathbb{T}}

\newcommand{\IP}{\mathbb{P}}

\newcommand{\drm}{\mathrm d}

\newcommand{\CD}{\mathcal D}

\newcommand{\CC}{\mathcal C}
\newcommand{\CH}{\mathcal H}

\newcommand{\CX}{\mathcal X}

\newcommand{\CB}{\mathcal B}

\newcommand{\SF}{\mathscr{F}}

\newcommand{\SL}{\mathscr{L}}

\renewcommand{\P}{\mathsf{P}}
\newcommand{\PT}{\widetilde{\P}}
\newcommand{\PI}{\mathsf{\Pi}}

\newcommand{\DC}{\mathsf{C}}

\newcommand{\DS}{\mathsf{S}}

\newcommand{\eps}{\varepsilon}

\definecolor{toon}{RGB}{71,145,67} 

\tikzset{baseline,odot/.style={shape=coordinate,circle,draw,inner sep=1pt},ori/.style={shape=coordinate,inner sep=0pt},cdot/.style={shape=coordinate,circle,fill,draw,inner sep=1pt},rdot/.style={shape=coordinate,rectangle,toon,draw,inner sep=1pt,outer sep=1.5pt,fill},idot/.style={shape=coordinate,rectangle,draw,inner sep=1pt,outer sep=1.5pt},pi/.style={decorate,decoration={zigzag,amplitude=2pt,pre length=4pt, post length=4pt,segment length=3pt}}}

\definecolor{ocre}{RGB}{64,123,121}
\definecolor{S}{rgb}{0.0,0.5,0.0}

\newcounter{item}
\numberwithin{item}{section}

\newtheorem{theorem}[item]{\sffamily Theorem}
\newtheorem{definition}[item]{\sffamily Definition}
\newtheorem{proposition}[item]{\sffamily Proposition}

\newtheorem{corollary}[item]{\sffamily Corollary}

\newtheorem*{theorem*}{\sffamily Theorem}
\newtheorem*{definition*}{\sffamily Definition}
\newtheorem*{proposition*}{\sffamily Proposition}
\newtheorem*{lemma*}{\sffamily Lemma}
\newtheorem*{corollary*}{\sffamily Corollary}

\usepackage[explicit]{titlesec}
\titleformat{\section}{\centering\Large\bfseries}{\thesection \ --}{0.7em}{\Large\bfseries #1}
\titleformat{\subsection}{\centering\large\bfseries}{\thesubsection \ --}{0.4em}{\large\bfseries #1}
\titleformat{\subsubsection}{\centering\bfseries}{\thesubsubsection \ --}{0.4em}{\bfseries #1}

\let\emph\relax
\DeclareTextFontCommand{\emph}{\bfseries\em}

\providecommand{\keywords}[1]
{
	{\footnotesize	
	\textbf{Keywords --} #1}
}

\title{\bfseries The infinitesimal generator of the Brox diffusion}
\author{Antoine MOUZARD\footnote{The author is supported by the Simons Collaboration on Wave Turbulence.}}
\date{}

\begin{document}

\maketitle
\abstract{We construct the infinitesimal generator of the Brox diffusion on a line with a periodic Brownian environment. This gives a new construction of the process and allows to solve the singular martingale problem. We prove that the associated semigroup is strong Feller with Gaussian lower and upper bounds. This also yields a construction of the Brox diffusion on a segment with periodic or Dirichlet boundary conditions. In this bounded space, we prove that there exists a unique measure and the existence of a spectral gap giving exponential ergodicity of the diffusion.}
\vspace{0.5cm}


\section*{Introduction}

The Brox diffusion is a random motion in random environment introduced by Brox in \cite{Brox}. It is the continuous analogue of Sinai's random walk \cite{Sinai82}, a random walk on the line $\IZ$ where the probability to move to the right at a point $x\in\IZ$ is given by $\omega_x$, and to the left by $1-\omega_x$, with $(\omega_x)_{x\in\IZ}$ independant and identically distributed random variables in $[0,1]$. The Brox diffusion is formaly given by
\begin{equation*}
\drm X_t=\nabla W(X_t)\drm t+\drm B_t
\end{equation*}
where $W$ is a two-sided standard Brownian motion independant of $B$ with $W(0)=0$. Its infinitesimal generator is formaly given by
\begin{equation*}
\SL=\frac{1}{2}\Delta+\nabla W\cdot\nabla
\end{equation*}
which is a singular stochastic operator due to the irregularity of $\nabla W$. To construct the diffusion, Brox rewrites the generator under the form
\begin{equation*}
\SL=\frac{1}{2e^{2W}}\frac{\drm}{\drm x}\left(\frac{1}{e^{-2W}}\frac{\drm}{\drm x}\right)
\end{equation*}
and consider the Itô-McKean construction using the self-similarity of the Brownian motion. In this work, we propose a new approach with the direct construction of the generator as a singular stochastic operator, in the case of a periodic Brownian motion $W$. While we choose to illustrate our method with this example, the method proposed here is general to solve SDEs with distributional drift and we choose to denote the derivative as $\frac{\drm}{\drm x}=\nabla$.

\medskip

The subject of SDEs with distributional drifts have recevied new attention with the recent development of rough paths and more generally regularity structures and paracontrolled calculus. It was originally studied by Mathieu \cite{Mathieu94} using Dirichlet forms methods where the drift corresponds to an irregular random media. Later, more general results were obtained by two different methods respectively by Bass and Chen \cite{BassChen01} and by Flandoli, Russo and Wolf \cite{FlandoliRussoWolf03,FlandoliRussoWolf04}, see also references therein. This work follows the second approach and consider the martingale problem associated to the SDE, taking its root in Stroock and Varadhan's book \cite{StroockVaradhan79}. For the SDE
\begin{equation*}
\drm X_t=b(X_t)\drm t+\sigma(X_t)\drm B_t
\end{equation*}
under suitable condition on the coefficient $b$ and $\sigma$, the Itô formula gives
\begin{equation*}
u(X_t)=u(X_0)+\int_0^t\big(\frac{1}{2}\sigma\sigma^*\Delta u+b\cdot\nabla u\big)(X_s)\drm s+\int_0^t\big(\sigma\nabla u\big)(X_s)\drm s
\end{equation*}
for smooth functions $u$. In particular, the process
\begin{equation}\label{star}\tag{$\star$}
u(X_t)-u(X_0)-\int_0^t\big(\frac{1}{2}\sigma\sigma^*\Delta u+b\cdot\nabla u\big)(X_s)\drm s
\end{equation}
is a martingale with respect to the filtration generated by $X$, and the differential operator
\begin{equation*}
\SL=\frac{1}{2}\sigma\sigma^*\Delta+b\cdot\nabla
\end{equation*}
is the infinitesimal generator of the process $X$. A martingale solution associated to an operator $(\SL,\CC)$ is a process $X$ such that (\ref{star}) defines a martingale for any function $u\in\CC$. This approach seems in particular interesting when $b$ is only a measurable function since the pointwise evaluation $b(X_t)$ does not make sense. As done in Stroock and Varadhan's book, this can be considered in the general case of a time dependent coefficient and the case of distributional time dependent drift was considered by Delarue and Diel \cite{DelarueDiel16} in the case $\sigma=1$, motivated by the recent progress on singular SPDEs. They considered the modified generator
\begin{equation*}
\partial_t+\frac{1}{2}\Delta+\nabla b_t\cdot\nabla
\end{equation*}
with $b_t$ an $\alpha$-Hölder function in space with $\alpha>\frac{1}{3}$ in one dimension. This was later extended by Cannizzaro and Chouk \cite{CC} on a torus $\IT^d$ for $d\le 3$ in the case of a random drift where they solve the same modified martingale problem. In particular, the backward Kolmogorov equation is a singular SPDE and they use the paracontrolled calculus and a renormalization procedure to go beyond the Young regime of \cite{DelarueDiel16}, relying on the recent progress on singular SPDEs. This is the reason why they consider the case of the torus, it is know to be easier to work in bounded space when dealing with white noise while it should be possible to adapt the techniques to the full space. Finally, Kremp and Perkowski \cite{KrempPerkowski22} generalized this approach to the case of a Lévy driving process. In this case, the Laplacian is replaced by the generator $\SL_\nu^\alpha$ of an $\alpha$-stable Lévy process with $\alpha\in(1,2]$ and they allow drift of spatial Hölder regularity $\frac{2-2\alpha}{3}$. It is expected that the regularity of the drift should depend on the regularizing properties of $\SL_\nu^\alpha$, they are exactly beyond the Young regime where a first order paracontrolled expansion is enough. As application, they solve the modified martingale problem associated to the Brox diffusion. They show that the condition
\begin{equation*}
u(X_t)-u(X_0)-\int_0^t(\partial_tu+\SL_\nu^\alpha u+\nabla W\cdot\nabla u)(X_s)\drm s
\end{equation*}
is a martingale for a sufficient large class of functions $u$ is a well-posed problem. In this work, we consider the martingale problem
\begin{equation*}
u(X_t)-u(X_0)-\int_0^t\Big(\frac{1}{2}\Delta u+\nabla W\cdot\nabla u\Big)(X_s)\drm s.
\end{equation*}
Note that the space $\CC$ of functions $u$ considered by Kremp and Perkowski is described using the paracontrolled calculus and does not contain smooth functions, the same will be true here and this is in large contrast with the classical case where one usually consider smooth compactly supported functions. A larger space $\CC$ imposes more condition hence making unicity more likely while one needs a small enough space $\CC$ to ensure existence. We identify here a suitable space $\CC$ such that there exists a unique martingale solution using the paracontrolled calculus.

\medskip

The study of singular stochastic operator of the form
\begin{equation*}
\Delta+b\cdot\nabla+a
\end{equation*}
for random fields $a,b$ was initiated by Allez and Chouk \cite{AllezChouk} following the recent development on singular SPDEs. They used the paracontrolled calculus and a renormalization procedure to construct the Anderson Hamiltonian
\begin{equation*}
\Delta+\xi
\end{equation*} 
where $\xi$ is the space white noise on the two dimensional torus $\IT^2$. This was later extended to various context by different authors, both with the paracontrolled calculus and regularity structures, see \cite{BDM,GUZ,Labbe,Mouzard,Ugurcan} and references therein. This is still a very active subject and going beyond the simple construction of the operator is a very hard problem. Another direction was followed by Morin and Mouzard \cite{MorinMouzard22} where they constructed the random magnetic Laplacian
\begin{equation*}
(i\nabla+A)^2=(i\partial_1+A_1)^2+(i\partial_2+A_2)^2
\end{equation*}
on a two dimensional torus $\IT^2$ with a potential $A:\IT^2\to\IR^2$ related to the white noise magnetic field $B=\xi$ via
\begin{equation*}
\xi=\nabla\times A=\partial_2A_1-\partial_1A_2.
\end{equation*}
Their construction relies on the heat paracontrolled calculus introduced in \cite{Mouzard} in the elliptic framework, see also \cite{BB1,BB3,BB2} for the parabolic calculus, which allows to deal with first order terms. With this work, we continue the study of singular stochastic operators. Indeed, the infinitesimal generator of the Brox diffusion
\begin{equation*}
\SL=\frac{1}{2}\Delta+\xi\cdot\nabla
\end{equation*}
on $\IT$ with $\xi=\nabla W$ the derivative of a Brownian motion on the torus falls in the range of the singular stochastic operators. We have $\xi\in\CC^{-\frac{1}{2}-\kappa}(\IT)$ for any $\kappa>0$ and this is just out of the range of the Young regime. We now explain the general construction of the domain of singular stochastic operators with the example of interest in this work, the infinitesimal generator of the Brox diffusion $\SL$. For smooth functions $u$, the operator is well-defined since the multiplication is not singular however one has
\begin{equation*}
\SL u=\frac{1}{2}\Delta u+\xi\cdot\nabla u\in\CH^{-\frac{1}{2}-\kappa}
\end{equation*}
which is a distribution. The idea is to consider in the domain rough functions $u$ with variations controlled by the noise such that $\frac{1}{2}\Delta u$ cancels the rough part of the product. Within the paracontrolled calculus, it is given by $\P_{\nabla u}\xi$ which has the same regularity as $\xi$ for $u\in\CH^1$. Motivated by a cancellation with the Laplacian, we can consider paracontrolled functions of the form
\begin{equation*}
u=\P_{\nabla u}X+u^\sharp
\end{equation*}
with $-\frac{1}{2}\Delta X=\xi$ and $u^\sharp$ a remainder smoother than $X$. Using the paracontrolled calculus, one can define the product $\xi\cdot\nabla u$ as long as one can define $\xi\cdot\nabla X$ which is a singular product. This is defined through a probabilist renormalization procedure, now usual in the resolution of singular SPDEs, and one works with the enhanced noise $\Xi=(\xi,\xi\cdot\nabla X)$. This yields that the operator $\SL$ is well-defined on such paracontrolled functions with values in $\CH^{-\kappa}$ for any $\kappa>0$. While this is not enough to define an unbounded operator in $L^2$, this gives a well-defined form
\begin{equation*}
(u,v)\mapsto\langle \SL u,v\rangle
\end{equation*}
for $u,v$ paracontrolled by $X$ since these are regular enough for the scalar product with the distribution in $\CH^{-\kappa}$ for $\kappa>0$ small to be well-defined, this is the form domain of $\SL$. Performing a higher order expansion yields for the operator $\SL$ the domain
\begin{equation*}
\{u\in L^2\ ;\ u-\P_{\nabla u}X_1-\P_{\nabla u}X_2\in\CH^2\}
\end{equation*}
with $X_1=X$ and $X_2=X_2(\Xi)$ a functionnal of the enhanced noise more regular than $X_1$. Once the domain is constructed, one can investigate any questions on the operator $(\SL,\CD_\Xi)$ and its spectral properties. This can also be implemented within regularity structures as done by Labbé \cite{Labbe} with different advantages and disadvantages

\medskip

In the first section, we construct the domain of the operator $\SL$ as paracontrolled functions. Using the heat paracontrolled calculus from \cite{Mouzard} and the $\Gamma$ map first introduced in \cite{GUZ}, the domain writes as $\CD_\Xi=\Gamma\CH^2$ where $\Xi$ corresponds to an enhancement of the noise $\xi$. We prove that $(\SL,\CD_\Xi)$ is closed, bounded from above and the limit of a suitable regularization $\SL_n$ in the resolvant sense. This is enough to define the semigroup associated to $-\SL$ which is conservative and strong Feller and we obtain Gaussian lower and upper bounds on its kernel. Conditionnaly to the environment $W$, this yields a Markov process $X$ generated by $\SL$ which corresponds to the Brox diffusion. We show that the martingale problem associated to $\SL$ with $\CC_\Xi=\Gamma C^2$ is well-posed, the solution being the law of the process $X$. We do not define the heat paracontrolled calculus and refer to \cite{Mouzard}, see however the Appendix for the needed continuity results on paraproducts and correctors. Up to this, our work is self-contained. Finally, this can be used to define the Brox diffusion on a segment with periodic or Dirichlet boundary condition where the generator has pure point spectrum with a spectral gap. We prove the existence of a unique invariant measure and exponential ergodicity of the diffusion follows from the spectral gap.

\medskip

As for the works \cite{CC,KrempPerkowski22}, the restriction of working with a noise in bounded spaces could be lifted at the price of working with weighted spaces. For the question of singular stochastic operators in unbounded space, we refer to Ugurcan's recent work \cite{Ugurcan} on the Anderson Hamiltonian. The method could also be adapted to deal with time dependent drift. The important feature of the process $W$ is its regularity. For example, one could deal with $W^H$ a fractionnal Brownian motion with index $H>\frac{1}{3}$ with the same methods and higher order expansion would be required to go beyond this regularity.

\medskip

For a smooth bounded potential $V:\IR^d\to\IR$, the Hilbert space $L^2(\IR^d,\drm\lambda)$ with $\lambda$ the Lebesgue measure is isometric to $L^2(\IT^2,e^{2V}\drm\lambda)$ and the generator
\begin{equation*}
Lu=\frac{1}{2}\Delta u+\nabla V\cdot\nabla u=\frac{1}{2}e^{-2V}\nabla\big(e^{2V}\nabla u\big)
\end{equation*}
is well-defined with domain $\CH^2$ and form domain $\CH^1$. We have
\begin{equation*}
\langle Lu,v\rangle_{L^2(\IR^d,e^{2V}\drm \lambda)}=-\frac{1}{2}\langle\nabla u,\nabla v\rangle_{L^2(\IR^d,e^{2V}\drm\lambda)}
\end{equation*}
thus $-L$ is a symmetric nonnegative operator in the Hilbert space $L^2(\IT^2,e^{2V}\drm\lambda)$. This can be used to show that $L$ is m-accretive thus generates a contraction semigroup $(e^{-tL})_{t\ge0}$ by the Hille-Yosida theorem with Gaussian lower and upper bounds depending on
\begin{equation*}
\delta(U)=\sup_{x\in\IR^d}V(x)-\inf_{x\in\IR^d}V(x),
\end{equation*} 
see section $4.3$ from Stroock's book \cite{stroock08}. This gives well-posedness of the martingale associated to the drift $\nabla V$ and this is the path we follow for the infinitesimal generator of the Brox diffusion $\SL$. If in addition the potential $V:\IR^d\to\IR$ is periodic, the differential operator
\begin{equation*}
\frac{1}{2}\Delta+\nabla V\cdot\nabla
\end{equation*}
can be considered on two different spaces, the periodic smooth functions $C_{\text{per}}^\infty(\IR^d)$ or the compactly supported smooth functions $C_0^\infty(\IR^d)$. This amounts to which Hilbert space we are working with, respectively $L^2(\IT^d)$ or $L^2(\IR^d)$, and this yields two differents semigroups. In the first case, the associated diffusion takes its values in $\IT^d$ while its takes its values in $\IR^d$ in the second case. The same holds for our construction of the Brox diffusion in a periodic environment and both approaches can be considered in the exact same way. We choose the second one which seems more natural, and this corresponds to the framework of Kremp and Perkowski \cite{KrempPerkowski22}. On the compact space $\IT^d$, one has stronger results on the spectral properties of the generator and we get a spectral gap for the generator and exponential ergodicity of the process to the unique invariant measure, this is the content of the last section.

\medskip

To construct $\SL$, we will need to solve $\Delta X=Y$ given $Y$, which does not have a unique solution. It is enough for our construction to consider a parametrix, one could consider for example the unique solution to $\Delta X=Y$ with $\int_{\IT}Y(x)\drm x=0$ or add a mass and consider $(\Delta+m)X=Y$ for any fixed $m>0$. To emphasize the generality of our methods, we choose to stick to the parametrix defined in \cite{Mouzard} and consider the somehow more complicated operator
\begin{equation*}
\Delta^{-1}:=\int_0^1e^{t\Delta}\drm t.
\end{equation*}
In particular, we have $\Delta\circ\Delta^{-1}=\text{Id}-e^{\Delta}$, an inverse up to the regularizing operator $e^{\Delta}$. We denote respectively as $\CC^\beta$ and $\CH^\beta$ the Besov-Hölder and Sobolev spaces on $\IR$ for $\beta\in\IR$.

\section{Domain of the generator}

The noise $\xi=\nabla W$ is given by
\begin{equation*}
\xi(x)=\sum_{k\in\IZ^*}\xi_ke^{ikx}
\end{equation*}
with $x\in\IR$, $(\xi_k)_{k\ge0}$ a family of independant and identically distributed random variable with law complex centered Gaussian with unit variance and $\xi_{-k}=\overline\xi_k$ for $k\ge1$ such that the noise is real. It satisfies the condition $\xi_0=0$ because the noise is the derivative of the periodic Brownian motion in the sense of distribution. Remark that a periodic Brownian motion is the nothing more than a Brownian bridge periodized over $\IR$. In particular, $\xi\in\CC^{-\frac{1}{2}-\kappa}$ for any $\kappa>0$ and we consider its regularization given by
\begin{equation*}
\xi_n(x)=\sum_{|k|\le n}\xi_ke^{ikx}
\end{equation*}
which converges to $\xi$ as $n$ goes to infinity in $\CC^{-\frac{1}{2}-\kappa}$. The regularized version of the SDE
\begin{equation*}
\drm X_t^{(n)}=\xi_n\big(X_t^{(n)}\big)\drm t+\drm B_t
\end{equation*}
admits a unique strong solution and the associated generator is
\begin{equation*}
\SL_n:=\frac{1}{2}\Delta+\xi_n\cdot\nabla
\end{equation*}
with domain $\CH^2$. In this section, we construct the operator $\SL$ with domain $\CD_\Xi$ as a limit of the operators $\SL_n$. Let $\alpha\in(1,\frac{3}{2})$ such that $\xi\in\CC^{\alpha-2}$ almost surely. We define the regularized enhanced noise as
\begin{equation*}
\Xi_n:=\big(\xi_n,\PI(\nabla X_n,\xi_n)\big)
\end{equation*}
with
\begin{equation*}
X_n:=-2\Delta^{-1}\xi_n.
\end{equation*}
In particular, it satisfies the equation
\begin{equation*}
-\frac{1}{2}\Delta X_n=\xi_n-e^\Delta\xi_n
\end{equation*}
and converges to $X=\Delta^{-1}\xi\in\CC^\alpha$ as $n$ goes to $\infty$. While the product $\nabla X\cdot\xi$ is not well-defined since $2\alpha-3<0$, the next proposition states that $\Xi_n$ converges to an enhancement of the noise $\Xi$ in its natural space 
\begin{equation*}
\CX^\alpha:=\CC^{\alpha-2}\times\CC^{2\alpha-3},
\end{equation*}
which is a product of distributions spaces.

\medskip

\begin{proposition}
There exist a distribution $\PI(\nabla X,\xi)\in\CC^{2\alpha-3}$ such that
\begin{equation*}
\lim_{n\to\infty}\big\|\big(\xi,\PI(\nabla X,\xi)\big)-\Xi_n\big\|_{\CX^\alpha}=0
\end{equation*}
and $\Xi:=\big(\xi,\PI(\nabla X,\xi)\big)\in\CX^\alpha$ is called the enhanced noise.
\end{proposition}

\medskip

\begin{proof}
This works as in the proof of lemma $5.1$ from \cite{KrempPerkowski22} and theorem $2.1$ from \cite{Mouzard}. We have $\xi\in\CC^{\alpha-2}$ almost surely thus the truncation in Fourier $\xi_n$ converges to $\xi$ in $\CC^{\alpha-2}$ as $n$ goes to infinity and we only have to show that $\big(\PI(\nabla X_n,\xi_n)\big)_{n\ge0}$ is a Cauchy sequence in $\CC^{2\alpha-3}$. Since the noise is periodic, it is enough to work on the torus $\IT$. The Besov spaces $\CB_{p,q}^\beta$ are defined by the norms
\begin{equation*}
\|u\|_{\CB_{p,q}^\beta}:=\Big(\sum_{j\ge-1}2^{\beta jq}\|\Delta_ju\|_{L^p(\IT)}^q\Big)^{\frac{1}{q}}
\end{equation*}
for $\beta\in\IR$ and $p,q\ge1$ with $\Delta_j$ the Littlewood-Paley projectors, see for example \cite{BahouriCheminDanchin11}. We have the particular cases $\CB_{2,2}^\beta=\CH^\beta$ and $\CB_{\infty,\infty}^\beta=\CC^\beta$ as well as the Besov embedding
\begin{equation*}
\CB_{p_1,q_1}^\beta\hookrightarrow\CB_{p_2,q_2}^{\beta-d(\frac{1}{p_1}-\frac{1}{p_2})}
\end{equation*}
for any $p_1\le p_2,q_1\le q_2$ and $\beta\in\IR$. We use 
\begin{equation*}
\CB_{2p,2p}^{\beta+\frac{1}{2p}}\hookrightarrow\CB_{\infty,\infty}^\beta=\CC^\beta
\end{equation*}
and bound the norm in $\CB_{2p,2p}^{\beta+\frac{1}{2p}}$ for a particular $p\ge1$ and $\beta=2\alpha-3$. For the Fourier resonant product
\begin{equation*}
Y_n:=\Pi(\nabla X_n,\xi_n)=\sum_{|\ell-\ell'|\le1}\Delta_\ell(\nabla X_n)\Delta_{\ell'}(\xi_n),
\end{equation*}
we have
\begin{align*}
\IE\big[\|Y_n-Y_m\|_{\CB_{2p,2p}^{\alpha+\frac{1}{2p}}}^{2p}\big]&=\sum_{j\ge-1}2^{(\alpha+\frac{1}{2p})j2p}\IE\big[\|\Delta_j(Y_n-Y_m)\|_{L^{2p}}^{2p}\big]\\
&=\sum_{j\ge-1}2^{(\alpha+\frac{1}{2p})j2p}\int_\IT\IE\big[|\Delta_j(Y_n-Y_m)(x)|^{2p}\big]\drm x\\
&\lesssim\sum_{j\ge-1}2^{(\alpha+\frac{1}{2p})j2p}\int_\IT\IE\big[|\Delta_j(Y_n-Y_m)(x)|^2\big]^p\drm x
\end{align*}
using Gaussian hypercontractivity. Since $Y_n-Y_m$ belongs to the second inhomogeneous Wiener chaos associated to $\xi$, we write
\begin{equation*}
\IE\big[|\Delta_j(Y_n-Y_m)(x)|^2\big]=\IE\big[|\Delta_j(Y_n-Y_m)(x)-c_{n,m}(x)|^2\big]+c_{n,m}(x)
\end{equation*}
where the first term belongs to the second homogeneous Wiener chaos and
\begin{equation*}
c_{n,m}(x)=\IE\big[\Delta_j(Y_n-Y_m)(x)\big].
\end{equation*}
The first term is arbitrary small uniformly with respect to $m\ge n\ge N$ for $N$ large enough as a Wick product, this is similar to the renormalization of the Anderson Hamiltonian. Since the law of the white noise is invariant by rotation, $c_{n,m}$ is a constant which usually correspond to a diverging constant. However in this case, we have $c_{n,m}(x)=0$ which can be seen with the computations
\begin{align*}
\IE\big[\nabla X_n(x)\cdot\xi_n(x)\big]&=\sum_{0<|k|,|k'|\le n}\frac{1}{ik}\IE\big[\xi_k\xi_{k'}\big]e^{i(k+k')x}\\
&=\sum_{0<|k|\le n}\frac{1}{ik}\\
&=0
\end{align*}
and using the symmetry of $\Delta_j$ which completes the proof for the Fourier resonant product. In this computation, we used $\IE[\xi_k^2]=0$ since $\xi_k$ is a complex centered Gaussian and $\IE[\xi_k\xi_{k'}]=0$ for $|k|\neq|k'|$. Note that $c_{n,m}(x)$ is also equal to $0$ for the heat resonant product $\PI$, the important properties is the symmetry of the truncation in frequencies, as in \cite{KrempPerkowski22}. We do not give more details in the computation to keep the presentation light and not introduce the details of the Paley-Littlewood projectors or the heat paracontrolled calculus. The proof for the heat resonant product can be directly adapted from \cite{Mouzard} theorem $2.1$.
\end{proof}

\begin{remark}
Our choice of regularization is important as one can see with the cancellation occuring in the previous proof. A different choice could be made at the price of substracting a diverging constant, as usually done with singular stochastic operators. For a general regularization $(\xi_\eps)_{\eps\ge0}$ of the noise, one should consider
\begin{equation*}
\PI(\nabla X,\xi):=\lim_{\eps\to0}\big(\PI(\nabla X_\eps,\xi_\eps)-c_\eps\big)
\end{equation*}
with $c_\eps=\IE\big[\PI(\nabla X_\eps,\xi_\eps\big]$ diverging like $|\log(\eps)|$ as $\eps$ goes to $0$. One has $c_\eps=0$ in the case of our symmetric regularization, as in the work of Kremp and Perkowski \cite{KrempPerkowski22}. One can see the possible logarithm divergence in the harmonic sum of the previous proof.
\end{remark}

\bigskip

We can now construct the domain $\CD_\Xi$ of $\SL$. It will consists of functions
\begin{equation*}
u=\PT_{\nabla u}X_1+\PT_{\nabla u}X_2+u^\sharp,
\end{equation*}
paracontrolled by $X_1=X\in\CC^\alpha$ and $X_2\in\CC^{2\alpha-1}$ solutions to
\begin{equation*}
-\frac{1}{2}\Delta X_1=\xi\quad\text{and}\quad-\frac{1}{2}\Delta X_2=\P_\xi\nabla X+\PI(\nabla X,\xi)
\end{equation*}
with $u^\sharp\in\CH^2$ a remainder. The paraproduct $\PT$ has the same properties as $\P$ and is intertwined via the relation
\begin{equation*}
\PT\circ\Delta^{-1}=\Delta^{-1}\circ\P,
\end{equation*}
see the appendix and \cite{Mouzard}. Since $\Delta\circ\Delta^{-1}=\text{Id}-e^{\Delta}$, one can first consider that $\Delta\circ\Delta^{-1}=\text{Id}$ when only interested in the regularity of the distributions and functions. As explained in the introduction, functions $u$ in the domain must have the correct form depending on the noise such that the term $\frac{1}{2}\Delta u$ cancel the rough part of the product $\xi\cdot\nabla u$. This choice gives
\begin{align*}
\SL u&=\frac{1}{2}\Delta\PT_{\nabla u} X_1+\frac{1}{2}\Delta\PT_{\nabla u}X_2+\frac{1}{2}\Delta u^\sharp+\xi\cdot\nabla u\\
&=-\P_{\nabla u}\xi-\P_{\nabla u}\big(\P_\xi\nabla X+\PI(\nabla X,\xi)\big)+\frac{1}{2}\Delta u^\sharp+\P_{\nabla u}\xi+\P_\xi\nabla u+\PI(\nabla u,\xi)\\
&=\frac{1}{2}\Delta u^\sharp-\P_{\nabla u}\P_\xi\nabla X-\P_{\nabla u}\PI(\nabla X,\xi)+\P_\xi\nabla u+\PI(\nabla u,\xi)
\end{align*}
thus $X_1$ cancel the paraproduct $\P_{\nabla u}\xi$. We use the paracontrolled toolkit from \cite{Mouzard}, see the appendix for the needed continuity results. Using the corrector $\DC_\nabla$, we have
\begin{align*}
\PI(\nabla u,\xi)&=\PI\big(\nabla\PT_{\nabla u}X_1,\xi\big)+\PI\big(\nabla\PT_{\nabla u}X_2,\xi\big)+\PI\big(\nabla u^\sharp,\xi\big)\\
&=\nabla u\cdot\PI(\nabla X_1,\xi)+\DC_\nabla(\nabla u,X_1,\xi)+\PI\big(\nabla\PT_{\nabla u}X_2,\xi\big)+\PI\big(\nabla u^\sharp,\xi\big)
\end{align*}
and the corrector $\DS$ gives
\begin{equation*}
\P_\xi\nabla u=\P_{\nabla u}\P_\xi\nabla X_1+\DS(\nabla u,X_1,\xi)+\P_\xi\nabla\PT_{\nabla u}X_2+\P_\xi u^\sharp
\end{equation*}
which gives the following definition.

\medskip

\begin{definition}
We define the domain of $\SL$ as
\begin{equation*}
\CD_\Xi:=\{u\in L^2\ ;\ u-\PT_{\nabla u}X_1-\PT_{\nabla u}X_2\in\CH^2\}.
\end{equation*}
For $u\in\CD_\Xi$, we have
\begin{equation*}
\SL u=\frac{1}{2}\Delta u^\sharp+\P_\xi u^\sharp+\PI(u^\sharp,\xi)+F_\Xi(u)
\end{equation*}
with
\begin{align*}
F_\Xi(u)&:=\P_{\PI(\nabla X_1,\xi)}\nabla u+\PI\big(\nabla u,\PI(\nabla X_1,\xi)\big)+\DC_\nabla(\nabla u,X_1,\xi)+\PI\big(\nabla\PT_{\nabla u}X_2,\xi\big)+\DS(\nabla u,X_1,\xi)+\P_\xi\nabla\PT_{\nabla u}X_2\\
&\quad-e^\Delta\big(\P_{\nabla u}X_1+\P_{\nabla u}X_2\big).
\end{align*}
\end{definition}

\medskip

At this point, it is not clear weither the domain is trivial or dense in $L^2$. Let $X_1^{(n)}$ and $X_2^{(n)}$ be the solutions of the same equations as $X_1$ and $X_2$ with $\Xi_n$ instead of $\Xi$. We introduce
\begin{equation*}
\Phi^{>N}(u):=u-\PT_{\nabla u}\big(X_1-X_1^{(N)}\big)-\PT_{\nabla u}\big(X_2-X_2^{(N)}\big)
\end{equation*}
where the notation is motivated by the fact that $X_i-X_i^{(N)}$ corresponds to the trunction to Fourier modes $|k|>N$. Since $X_1^{(N)}$ and $X_2^{(N)}$ are smooth, we have
\begin{equation*}
\CD_\Xi=(\Phi^{>N})^{-1}(\CH^2)
\end{equation*}
for any $N\ge0$, the case $N=0$ corresponding to the definition. Since
\begin{equation*}
\lim_{N\to\infty}\|X_1-X_1^{(N)}\|_{\CC^\alpha}+\|X_2-X_2^{(N)}\|_{\CC^{2\alpha-1}}=0,
\end{equation*}
the map $\Phi^{>N}:\CH^\beta\to\CH^\beta$ is invertible for $N$ large enough and any $\beta\in[0,\alpha)$ as a small perturbation of the identity. Indeed, we have
\begin{align*}
\|u-\Phi^{>N}(u)\|_{\CH^\beta}&\lesssim\|\nabla u\|_{L^2}\big(\|X_1-X_1^{(N)}\|_{\CC^\alpha}+\|X_2-X_2^{(N)}\|_{\CC^\alpha}\big)\\
&\lesssim\|u\|_{\CH^\beta}\big(\|X_1-X_1^{(N)}\|_{\CC^\alpha}+\|X_2-X_2^{(N)}\|_{\CC^{2\alpha-1}}\big)
\end{align*}
for $\beta\in(1,\alpha)$ and
\begin{align*}
\|u-\Phi^{>N}(u)\|_{\CH^\beta}&\lesssim\|\nabla u\|_{\CH^{\beta-1}}\big(\|X_1-X_1^{(N)}\|_{\CC^\alpha}+\|X_2-X_2^{(N)}\|_{\CC^\alpha}\big)\\
&\lesssim\|u\|_{\CH^\beta}\big(\|X_1-X_1^{(N)}\|_{\CC^\alpha}+\|X_2-X_2^{(N)}\|_{\CC^{2\alpha-1}}\big)
\end{align*}
for $\beta\in[0,1]$ using that $\alpha>1$. For $N\ge N_\Xi$ with $N_\Xi$ depending only on the norm of the enhanced noise, the map $\Phi^{>N}:\CH^\beta\to\CH^\beta$ is invertible for $\beta\in[0,\alpha)$ as a small perturbation of the identity, we denote as $\Gamma^{>N}$ its inverse. It is defined by the implicit equation
\begin{equation*}
\Gamma^{>N}u^\sharp=\PT_{\nabla\Gamma^{>N}u^\sharp}\big(X_1-X_1^{(N)}\big)+\PT_{\nabla\Gamma^{>N}u^\sharp}\big(X_2-X_2^{(N)}\big)+u^\sharp
\end{equation*}
for any $u^\sharp\in\CH^\beta$. This implies that $\CD_\Xi=\Gamma^{>N}\CH^2$ is not trivial, in particular dense with the following proposition. In the following, we will denote as $\Phi=\Phi^{>N_\Xi}$ and $\Gamma=\Gamma^{>N_\Xi}$ in order to simplify the notation but the reader should keep in mind that it depends on the size of the enhanced noise $\|\Xi\|_{\CX^\alpha}$, a random quantity. The idea of introducing this map $\Gamma$ was introduced by Gubinelli, Ugurcan and Zachhuber in \cite{GUZ}. Finally, we have a natural regularization for functions in the domain. Indeed, introduce
\begin{equation*}
\Phi_n(u):=u-\PT_{\nabla u}\big(X_1^{(n)}-X_1^{(N_\Xi)}\big)-\PT_{\nabla u}\big(X_2^{(n)}-X_2^{(N_\Xi)}\big)
\end{equation*}
and its inverse $\Gamma_n$. Then we have that $\Gamma_nu^\sharp\in\CH^\beta$ for $u^\sharp\in\CH^\beta$ for any $\beta\in\IR$ and that $\Phi_n$ and $\Gamma_n$ respectively converge to $\Phi$ and $\Gamma$ as $n$ goes to infinity, see the Appendix for the precise results. This indeed gives a regularisation of the domain since
\begin{equation*}
\lim_{n\to\infty}\|\Gamma u^\sharp-\Gamma_nu^\sharp\|_{\CH^\beta}=0
\end{equation*}
for any $u^\sharp\in\CH^2$ and $\beta\in[0,\alpha)$.

\medskip

\begin{proposition}
The domain $\CD_\Xi$ is dense in $\CH^\beta$ for any $\beta\in[0,\alpha)$.
\end{proposition}

\medskip

\begin{proof}
Let $f\in\CH^2$ and consider $g_n:=\Phi_n(f)\in\CH^2$. We have
\begin{equation*}
\|f-\Gamma(g_n)\|_{\CH^\beta}=\|\Gamma_n(g_n)-\Gamma(g_n)\|_{\CH^\beta}\le\|\Gamma_n-\Gamma\|_{\CH^\beta\to\CH^\beta}\|g_n\|_{\CH^\beta}
\end{equation*}
for any $\beta\in[0,\alpha)$. Using that $\Gamma_n$ converges to $\Gamma$ in norm as bounded operators on $\CH^\beta$ and that $\Phi_n$ is bounded on $\CH^\beta$ uniformly with respect to $n$, we get
\begin{equation*}
\lim_{n\to\infty}\|f-\Gamma(g_n)\|_{\CH^\beta}=0.
\end{equation*}
Since $g_n\in\CH^2$, we have $\Gamma(g_n)\in\CD_\Xi$ and the proof is complete using that $\CH^2$ is dense in $\CH^\beta$.
\end{proof}

\begin{remark}
While smooth functions are usally in the domain of a differential operator, this is not the case for singular stochastic operator. This will be important when considering the martingale problem associated to $\SL$. The precise Hölder regularity of functions in the domain is computed below. However, in this case, we have $1=\Gamma 1$ thus constant functions are in the domain. Having in mind the martingale problem, this is natural since constant processes are indeed martingales.
\end{remark}

\bigskip

We end this section with the computation of the Hölder regularity of functions in the domain, which is more natural for probabilist. In particular, functions of the domain have continuous derivative.

\medskip

\begin{proposition}
We have
\begin{equation*}
\CD_\Xi\subset\CC^{\frac{3}{2}-\kappa}
\end{equation*}
for any $\kappa>0$.
\end{proposition}

\medskip

\begin{proof}
In one dimension, we have
\begin{equation*}
\CH^\beta\hookrightarrow\CC^{\beta-\frac{1}{2}}
\end{equation*}
for any $\beta\in\IR$. Since $\CD_\Xi=\Phi^{-1}(\CH^2)$ and $\Phi:\CH^\beta\to\CH^\beta$ is invertible for any $\beta\in[0,\alpha)$, we get
\begin{equation*}
\nabla\CD_\Xi\subset\CC^{\alpha-\frac{3}{2}}\subset\CC^{-\kappa}
\end{equation*}
for any $\kappa>0$. Since $X_1,X_2\in\CC^\alpha$ and $\CH^2\hookrightarrow\CC^\alpha$, we get the result.
\end{proof}

\section{Properties of the generator}

With the parametrization of the domain, we have the natural norm
\begin{equation*}
\|u\|_{\CD_\Xi}^2=\|u\|_{L^2}^2+\|\Phi(u)\|_{\CH^2}^2
\end{equation*}
which is in fact equivalent to the graph norm
\begin{equation*}
\|u\|_\SL^2=\|u\|_{L^2}^2+\|\SL u\|_{L^2}^2.
\end{equation*}
This guarantees that $(\SL,\CD_\Xi)$ is a closed operator.

\medskip

\begin{proposition}
The operator $(\SL,\CD_\Xi)$ is closed in $L^2$.
\end{proposition}

\medskip

\begin{proof}
Let $(u_n)_{n\ge1}\subset\CD_\Xi$ such that
\begin{equation*}
\lim_{n\to\infty}\|u_n-u\|_{L^2}^2+\|\SL u_n-v\|_{L^2}^2=0
\end{equation*}
for $u,v\in L^2$. Since $u_n\in\CD_\Xi$, we have
\begin{equation*}
\|\Phi(u_n)\|_{\CH^2}\lesssim\|u_n\|_{L^2}+\|\SL u_n\|_{L^2},
\end{equation*}
thus $(\Phi(u_n)\big)_{n\ge0}$ is a Cauchy sequence in $\CH^2$ and converges to a function $u^\sharp\in\CH^2$. Since $\Phi$ is continuous on $L^2$, we have $\Phi(u)=u^\sharp$ and $u\in\CD_\Xi$. Finally, we have
\begin{align*}
\|\SL u-v\|_{L^2}&\le\|\SL u-\SL u_n\|_{L^2}+\|\SL u_n-v\|_{L^2}\\
&\lesssim\|u-u_n\|_{L^2}+\|u_n^\sharp-u^\sharp\|_{\CH^2}+\|\SL u_n-v\|_{L^2}
\end{align*}
thus $\SL u=v$ and the operator is closed.
\end{proof}

We constructed the domain of the operator, that is a space $\CD_\Xi\subset L^2$ such that $\SL$ sends $\CD_\Xi$ in $L^2$. As explained in the introduction, we also get the form domain, a space $\CD(\sqrt{L})$ such that $\CD_\Xi\subset\CD(\sqrt{L})\subset L^2$ and $\langle\SL u,u\rangle$ is well-defined for any $u\in\CD(\sqrt{\SL})$. The form domain of $\SL$ is
\begin{equation*}
\CD^\beta\big(\sqrt{\SL}\big):=\{u\in L^2\ ;\ u-\PT_{\nabla u}X_1\in\CH^\beta\}
\end{equation*}
for any $\beta>\frac{3}{2}$. As for $\CD_\Xi$, the form domain is parametrized by $\CH^\beta$ with
\begin{equation*}
\CD^\beta\big(\sqrt{\SL}\big)=\Gamma\CH^\beta.
\end{equation*}
We prove that the operator $\SL_n\Gamma_n$ converges to $\SL\Gamma$ in norm and a convergence of form.

\medskip

\begin{proposition}
For any $u^\sharp\in\CH^2$, we have
\begin{equation*}
\|\SL\Gamma u^\sharp-\SL_n\Gamma_nu^\sharp\|_{L^2}\lesssim\|u^\sharp\|_{\CH^2}\|\Xi-\Xi_n\|_{\CX^\alpha}.
\end{equation*}
Moreover, for any $u^\sharp\in\CH^\beta$ and $\beta>\frac{3}{2}$, we have
\begin{equation*}
\lim_{n\to\infty}\big\langle\SL_n\Gamma_nu^\sharp,\Gamma_nu^\sharp\big\rangle=\big\langle\SL\Gamma u^\sharp,\Gamma u^\sharp\big\rangle.
\end{equation*}
\end{proposition}

\medskip

\begin{proof}
Let $u^\sharp\in\CH^2$ and consider $u:=\Gamma u^\sharp$ and $u_n:=\Gamma_nu^\sharp$. We have
\begin{equation*}
\SL u=\frac{1}{2}\Delta u^\sharp+\P_\xi\nabla u^\sharp+\PI\big(\nabla u^\sharp,\xi\big)+F_\Xi(u)
\end{equation*}
with
\begin{align*}
F_\Xi(u)&=\P_{\PI(\nabla X_1,\xi)}\nabla u+\PI\big(\nabla u,\PI(\nabla X_1,\xi)\big)+\DC(\nabla u,X_1,\xi)+\PI\big(\nabla\PT_{\nabla u}X_2,\xi\big)+\DS(\nabla u,X_1,\xi)+\P_\xi\nabla\PT_{\nabla u}X_2\\
&\quad-e^\Delta\big(\P_{\nabla u}X_1+\P_{\nabla u}X_2\big).
\end{align*}
Thus
\begin{equation*}
\SL u-\SL_nu_n=F_\Xi(u)-F_{\Xi_n}(u_n)
\end{equation*}
and we only have to control the nonlinear term $F_\Xi$, which is granted using the continuity results of the paracontrolled theory and
\begin{equation*}
\lim_{n\to\infty}\|u-u_n\|_{\CH^\gamma}+\|\Xi-\Xi_n\|_{\CX^\alpha}=0
\end{equation*}
for any $\gamma\in[0,\alpha)$. For the second part, let $u^\sharp\in\CH^\beta$ and consider again $u=\Gamma u^\sharp$ and $u_n=\Gamma_nu^\sharp$. We have
\begin{equation*}
\SL u=\frac{1}{2}\Delta u^\sharp+\P_\xi\nabla u+\nabla u\PI(\nabla X_1,\xi)+\DC_\nabla(\nabla u,X_1,\xi)+\PI(\nabla u^\sharp,\xi)-e^{\Delta}\P_{\nabla u}X_1
\end{equation*}
which is a well-defined distribution in $\CH^{\beta-2\wedge2\alpha-3}$. Since $u\in\CH^\alpha$ with $\alpha=\frac{3}{2}-\kappa$ for $\kappa$ arbitrarty small, the scalar product $\langle\SL u,u\rangle$ is indeed well-defined. The convergence of $\Xi_n$ to $\Xi$ in $\CX^\alpha$ guarantees that $\SL_nu_n$ converges to $\SL u$ and this completes the proof.
\end{proof}

In order to define the semigroup associated to $-\SL$, we want to show that the form is bounded from below, that is the existence of a possibly random constant $c=c(\Xi)>0$ such that
\begin{equation*}
\forall u\in\CD(\sqrt{\SL}),\quad\langle-\SL u,u\rangle\ge-c\|u\|_{L^2}^2.
\end{equation*}
The method used in \cite{GUZ,MorinMouzard22,Mouzard} for the Anderson Hamiltonian and the random magnetic Laplacian can not be used here since the remainder of the paracontrolled expansion giving the form domain needs to be of Sobolev regularity strictly greater than $1$. Also, $\SL$ is not symmetric nor selfadjoint in $L^2(\IR,\drm\lambda)$. However, even with regularized drift $\nabla W_n$, it is useful to rewrite the operator as
\begin{equation*}
-\SL_nv=-\frac{1}{2}e^{-2W_n}\nabla\left(e^{2W_n}\nabla v\right)
\end{equation*}
which is symmetric in the Hilbert space $L^2(\IR,e^{2W_n}\drm\lambda)$. Moreover, the associated form is nonnegative with
\begin{equation*}
\langle-\SL_nv,v\rangle_{L^2(\IR,e^{2W_n}\drm\lambda)}=\langle\nabla v,\nabla v\rangle_{L^2(\IR,e^{2W_n}\drm\lambda)}\ge0
\end{equation*}
for any $v\in\CH^1$. Using the convergence of $W_n$ to $W$ and the previous proposition, this gives that $-\SL$ is nonnegative in the Hilbert space $L^2(\IR,e^{2W}\drm\lambda)$, the regularity of $W$ being enough for the scalar product to make sense. Since $W$ is a continuous bounded function, the spaces $L^2(\IR,\drm\lambda)$ and $L^2(\IR,e^{2W_n}\drm\lambda)$ are isometrics and we can state continuity results and work with functions and norm considered in either spaces.

\medskip

\begin{corollary}
For any $\beta>\frac{3}{2}$ and $u\in\CD^\beta(\sqrt{\SL})$, we have
\begin{equation*}
\langle-\SL u,u\rangle_{L^2(\IR,e^{2W}\drm\lambda)}=\langle-e^{2W}\SL u,u\rangle\ge0,
\end{equation*}
that is $-\SL$ is nonnegative in the Hilbert space $L^2(\IR,e^{2W}\drm\lambda)$.
\end{corollary}

\medskip

The notation $\CD^\beta(\sqrt{-\SL})$ for the form domain would be more accurate since $-\SL\ge0$, we keep the previous one for simplicity. Applying the Babuška–Lax–Milgram yields the following theorem.

\medskip

\begin{proposition}
For any $c>0$, the operators $\SL-c$ and $\SL_n-c$ are invertibles and
\begin{align*}
(\SL-c)^{-1}&:L^2\to\CD_\Xi,\\
(\SL_n-c)^{-1}&:L^2\to\CH^2
\end{align*}
are bounded.
\end{proposition}

\medskip

\begin{proof}
Following \cite{GUZ,MorinMouzard22,Mouzard}, we want to use the theorem of Babuška-Lax-Milgram from \cite{Babuska}, a generalization of the Lax-Milgram theorem. In this proof, we write $\|\cdot\|_{L^2}$ for the norm in $L^2(\IT,e^{2W}\drm\lambda)$ since we are working in this Hilbert space. For any $c>0$, we have
\begin{equation*}
c\|u\|_{L^2}^2<\big\langle(-\SL+c)u,u\big\rangle_{L^2(\IT,e^{2W}\drm\lambda)}
\end{equation*}
for any $u\in\CD_\Xi$ using the previous corollary. Consider the norm
\begin{equation*}
\|u\|_{\CD_\Xi}^2=\|u\|_{L^2}^2+\|u^\sharp\|_{\CH^2}^2
\end{equation*}
with $u=\Gamma u^\sharp$ on the space $\CD_\Xi$. This yields that $-\SL+c$ is a weakly coercive operator, that is
\begin{equation*}
c\|u\|_{\CD_\Xi}\le\|(-\SL+c)u\|_{L^2}=\sup_{\|v\|_{L^2}=1}\big\langle(-\SL+c)u,v\big\rangle_{L^2(\IT,e^{2W}\drm\lambda)}
\end{equation*}
for any $u\in\CD_\Xi$. Moreover, the bilinear map
\begin{equation*}
\left.\begin{array}{cccc}
B:&\CD_\Xi\times L^2&\to&\IR\\
&(u,v)&\mapsto&\big\langle(-\SL+c)u,v\big\rangle_{L^2(\IT,e^{2W}\drm\lambda)}
\end{array}\right.
\end{equation*}
is continuous since
\begin{equation*}
\quad|B(u,v)|\le\|(-\SL+c)u\|_{L^2}\|v\|_{L^2}\lesssim_\Xi\|u\|_{\CD_\Xi}\|v\|_{L^2}
\end{equation*}
for $u\in\CD_\Xi$ and $v\in L^2$. The last condition we need is that for any $v\in L^2\backslash\{0\}$, we have
\begin{equation*}
\sup_{\|u\|_{\CD_\Xi}=1}|B(u,v)|>0.
\end{equation*}
Let assume that there exists $v\in L^2$ such that $B(u,v)=0$ for all $u\in\CD_\Xi$. Then
\begin{equation*}
\forall u\in\CD_\Xi,\quad\langle u,v\rangle_{\CD_\Xi,\CD_\Xi^*}=0.
\end{equation*}
hence $v=0$ as an element of $\CD_\Xi^*$. By density of $\CD_\Xi$ in $L^2(\IT,e^{2W}\drm\lambda)$, this implies $v=0$ in $L^2$ hence the needed property. By the theorem of Babuška-Lax-Milgram, for any $f\in L^2$ there exists a unique $u\in\CD_\Xi$ such that
\begin{equation*}
\forall v\in L^2,\quad B(u,v)=\langle f,v\rangle_{L^2(\IT,e^{2W}\drm\lambda)}.
\end{equation*}
Moreover, we have $\|u\|_{\CD_\Xi}\lesssim_\Xi\|f\|_{L^2}$ hence the result for $(\SL-c)^{-1}$. The same argument works for $\SL_n-c$ since $-\SL_n$ is also nonnegative in $L^2(\IR,e^{2W}\drm\lambda)$ for any $n\ge1$ with domain $\CH^2$.
\end{proof}

We have now proven that $-\SL$ is a closed $m$-accretive operator hence generates a contraction semigroup $(e^{-t\SL})_{t\ge0}$ by the Hille-Yosida theorem. We have also shown that the operator $\SL$ is natural as the limit of $\SL_n$ in some sense. However, the domain $\CD_\Xi$ of $\SL$ is disjoint of $\CH^2$, the domain of $\SL_n$, thus one can not compare $\SL u$ and $\SL_ nu$. The correct limits to consider is $\SL_n\Gamma_n$ to $\SL\Gamma$ as done before or the convergence of the resolvants, given by the following proposition.

\medskip

\begin{proposition}
For any $\beta\in[0,\alpha)$ and $c>0$, we have
\begin{equation*}
\|(\SL+c)^{-1}-(\SL_n+c)^{-1}\|_{\CH^\beta\to\CH^\beta}\lesssim\|\Xi-\Xi_n\|_{\CX^\alpha}.
\end{equation*}
We say that $\SL_n$ converges to $\SL$ in the resolvant sense.
\end{proposition}

\medskip

\begin{proof}
Let $v\in L^2$ and $c>0$. Since $\SL+c:\CD_\Xi\to L^2$ is invertible, there exist $u\in\CD_\Xi$ such that
\begin{equation*}
(\SL+c)u=v.
\end{equation*}
Let $u^\sharp=\Phi(u)$ and $u_n=\Gamma_nu^\sharp$. We have
\begin{align*}
\big\|(\SL+c)^{-1}v-(\SL_n+c)^{-1}v\big\|_{\CH^\beta}&=\big\|u-(\SL_n+c)^{-1}(\SL+c)u\big\|_{\CH^\beta}\\
&\le\|u-u_n\|_{\CH^\beta}+\big\|u_n-(\SL_n)^{-1}(\SL+c)u\big\|_{\CH^\beta}\\\
&\le\|u-u_n\|_{\CH^\beta}+\big\|(\SL_n+c)^{-1}\big((\SL_n+c)u_n-(\SL+c)u\big)\big\|_{\CH^\beta}
\end{align*}
which completes the proof using the convergence in norm of $\SL_n\Gamma_n$ to $\SL\Gamma$ and that $(\SL_n+c)^{-1}$ is uniformly bounded with respect to $n$ from $L^2$ to $\CH^\beta$.
\end{proof}

For numerical function, we have
\begin{equation*}
e^{-t}=\lim_{N\to\infty}\Big(1+\frac{t}{N}\Big)^{-N}=e^{tc}\lim_{n\to\infty}\Big(1+\frac{t}{N}(1+c)\Big)^{-N}
\end{equation*}
for any $t,c\in\IR$. Thus we have
\begin{equation*}
e^{t\SL}=e^{tc}\lim_{N\to\infty}\Big(1+\frac{t}{N}(\SL+c)\Big)^{-N}
\end{equation*}
as an operator from $L^2$ to $\CD_\Xi$ using the convergence of the resolvant, this can also be seen as a definition of the semigroup generated by $\SL$. We denote respectively as $p_t(x,y)$ and $p_t^{(n)}(x,y)$ the kernel of $e^{t\SL}$ and $e^{t\SL_n}$ for $x,y\in\IR$ and $t>0$. For any $n>0$, the regularized operator
\begin{equation*}
\SL_n=\frac{1}{2}\Delta+\xi_n\cdot\nabla
\end{equation*}
is a conservative perturbation of the Laplacian since $\xi_n=\nabla W_n$ with $W_n$ the truncation in Fourier of the periodic Brownian motion $W$. Following section $4.3$ from Stroock's book \cite{stroock08}, we get uniform lower and upper bounds on the Gaussian bounds on the associated semigroup. This yields the following theorem. Note that since $\SL1=0$, the semigroup is conservative with $\int_{\IR}p_t(x,y)\drm y=1$ for any $t>0$ and $x\in\IR$.

\medskip

\begin{theorem}
There exists random constants $c=c(\xi),m=m(\xi)>0$ such that
\begin{equation*}
\frac{1}{c\sqrt{t}}e^{-c\frac{|x-y|^2}{t}}\le p_t^{(n)}(x,y),p_t(x,y)\le\frac{c}{\sqrt{t}}e^{-\frac{|x-y|^2}{ct}}
\end{equation*}
for any $x,y\in\IR$, $t>0$ and $n\ge0$.
\end{theorem}

\medskip

\begin{proof}
For the regularized operator $\SL_n$, the result follows from Stroock's book \cite{stroock08} with constants $c_n,m_n$ depending only on
\begin{equation*}
\delta(W_n):=\sup_{x\in\IT}W_n(x)-\inf_{x\in\IT}W_n(x).
\end{equation*}
In particular, we have
\begin{equation*}
\delta(W_n)\le2\sup_{x\in\IT}\Big|W_n(x)-\int_{\IT}W_n(x)\drm x\Big|
\end{equation*}
hence $\delta(W_n)$ is determined by $\xi_n$. For any $\kappa>0$, $W\in\CC^{\frac{1}{2}-\kappa}$ almost surely thus $W_n$ converges to $W$ in $\CC^{\frac{1}{2}-\kappa}$ almost surely. In particular, this implies
\begin{equation*}
\sup_{n\ge0}\|W_n\|_{\CC^{\frac{1}{2}-\kappa}}\le2\|W_n\|_{\CC^{\frac{1}{2}-\kappa}}<\infty
\end{equation*}
hence there exists constants $c_\xi,m_\xi>0$ depending only on $\xi$ such that
\begin{equation*}
\frac{1}{c\sqrt{t}}e^{-c\frac{|x-y|^2}{t}}\le p_t^{(n)}(x,y)\le\frac{c}{\sqrt{t}}e^{-\frac{|x-y|^2}{ct}}
\end{equation*}
for any $x,y\in\IT$, $t>0$ and $n\ge0$. The results follows using the convergence the semigroups of $\SL_n$ to the one of $\SL$.
\end{proof}

\section{The martingale problem}

For random dynamics in random environment, one has two different sources of randomness. Here, this corresponds to the random variable $\xi$ that gives the environment and the Brownian motion $B$ independent of $\xi$. One can consider the random variable $X$ where both the environment and the driving path are random, this is called the annealed law, or conditionned to the environment corresponding to the quenched law. The construction of the random stochastic operator $\SL$ depends on the environment and, conditionnaly to $\xi$, the semigroup $e^{t\SL}$ generates a Markov process $(X_t)_{t\ge0}$ where the finite dimensional distributions
\begin{equation*}
\IP(X_{t_1}\in A_1,\ldots,X_{t_n}\in A_n)=\int_{A_1\times\ldots\times A_n}\left(\prod_{i=1}^{n-1}p_{t_{i+1}-t_i}(x_i,x_{i+1})\right)\drm x_1\ldots\drm x_n
\end{equation*}
for any $0<t_1\le\ldots\le t_n$ and $A_1,\ldots,A_n$ measurable sets of $\IR$. This yields a measure on $\SF(\IR^+,\IR)$ and the upper Gaussian bounds together with Kolmogorov's theorem ensures that it has the same Hölder regularity as the Brownian motion.

\medskip

\begin{proposition}
The path of the process $X$ belongs almost surely to $\CC^{\frac{1}{2}-\kappa}$ for any $\kappa>0$.
\end{proposition}

\medskip

The lower Gaussian bounds guarantee that the diffusion spreads to the whole space with positive probability at any time $t>0$. Moreover, since the semigroup $e^{-t\SL}$ is continuous from $L^2$ to the domain $\CD_\Xi$ for any $t>0$, it is strong Feller in the sense that is sends measurables bounded functions to continuous functions. In this section, we show that this process is the Brox diffusion and its law is caracterized as the unique solution to the martingale associated to $\SL$.

\medskip

\begin{definition}
A law $\IQ$ on $C\big(\IR^+,\IR\big)$ is called a solution to the martingale problem associated to $(\SL,\CC)$ with initial data $x\in\IR$ if for any process $X$ of law $\IQ$ and function $u\in\CC$, the process
\begin{equation*}
u(X_t)-u(X_0)-\int_0^t(\SL u)(X_s)\drm s
\end{equation*}
is a martingale with respect to the filtration generated by $X$ and $X_0=x$ almost surely.
\end{definition}

\medskip

Usually, the choice of the space $\CC\subset\CD(\SL)$ is not important and one considers the space of smooth and compactly supported functions which are always in the domain of a differential operator. For the Brox diffusion, smooth functions are not in the domain of the generator and the choice of a suitable space is one of the achievements of this work. As for the Laplacian with $\CH^2$, the domain is too large in order to guarantee existence and one has to choose a core. We consider
\begin{equation*}
\CC_\Xi:=\Gamma\CC^2\subset\CD_\Xi
\end{equation*}
and prove that the martingale problem associated to $(\SL,\CC_\Xi)$ is well-posed.

\medskip

\begin{theorem}
There exists a unique solution to the martingale problem associated to $(\SL,\CC_\Xi)$. Moreoever, the solution is the weak limit of the solutions $X^{(n)}$ of the regularized equation
\begin{equation*}
\drm X_t^{(n)}=\xi_n\big(X_t^{(n)}\big)\drm t+\drm B_t.
\end{equation*}
\end{theorem}

\medskip

\begin{proof}
Using the upper Gaussian bounds on the semigroup $e^{t\SL_n}$ uniform with respect to $n$, there exists a constant $c=c(\xi)>0$ independent of $n\ge0$ such that
\begin{equation*}
\IE\big[|X_t^{(n)}-X_s^{(n)}|^p\big]\le c|t-s|^{\frac{p}{2}}
\end{equation*}
for any $0\le s\le t$. Kolmogorov's tightness criterion implies that $X^{(n)}$ is tight on $C\big(\IR^+,\IR\big)$ and, up to extraction, we suppose that the law of $X^{(n)}$ converges to a law $\IQ$ on $C\big(\IR^+,\IR\big)$. Moreoever, the convergence of the semigroup $(e^{t\SL_n})_{t\ge0}$ to the semigroup generated by $\SL$ implies the convergence of the finite dimensional distributions of $X^{(n)}$ to $X$ thus $\IQ$ is the law of $X$, the Markov process with generator $\SL$. We now prove that $X$ is a solution to the martingale problem associated to $(\SL,\CD_\Xi)$. Let $u^\sharp\in\CC^2$ and consider $u_n:=\Gamma_nu^\sharp\in\CC^2$ and
\begin{equation*}
M_t^n(Y):=u_n(Y_t)-u_n(Y_0)-\int_0^t(\SL_nu_n)(Y_r)\drm r
\end{equation*}
for any path $Y\in C\big([0,T],\IR\big)$. Then $\big(M_t^n(X^{(n)})\big)_{0\le t\le T}$ is a martingale with respect to the filtration generated by $X^{(n)}$. Thus for any $0\le s\le t$ and $F:C([0,s],\IR)\to\IR$ continuous bounded function, we have
\begin{equation*}
\IE\Big[\big(M_t^n(X^{(n)})-M_s^n(X^{(n)})\big)F\big((X_r^{(n)})_{0\le r\le s}\big)\Big]=0.
\end{equation*}
Our goal is to prove that
\begin{equation*}
\IE\Big[\big(M_t(X)-M_s(X)\big)F\big((X_r)_{0\le r\le s}\big)\Big]=0
\end{equation*}
for $u=\Gamma u^\sharp\in\CD_\Xi$ to conclude that the process $\big(M_t(u)\big)_{0\le t\le T}$ is a martingale with respect to the filtration generated by $X$. Using Itô formula, we have
\begin{equation*}
\IE\big[|M_t^n(X^{(n)})|^2\big]\le T\|\nabla u_n\|_{L^\infty}^2\le 2T\|\nabla u\|_{L^\infty}
\end{equation*}
for $n$ large enough using that $u_n$ converges to $u$ in $\CC^\alpha$ with $\alpha>1$. Since $F$ is bounded, this implies that the family $\big(M_t^n(X^{(n)})-M_s^n(X^{(n)})\big)F\big((X_r^{(n)})_{0\le r\le s}\big)_{n\ge0}$ is uniformly integrable hence we can invert integration and limit to get
\begin{equation*}
\lim_{n\to\infty}\IE\Big[\big(M_t^n(X^{(n)})-M_s^n(X^{(n)})\big)F\big((X_r^{(n)})_{0\le r\le s}\big)\Big]=\IE\Big[\big(M_t(X)-M_s(X)\big)F\big((X_r)_{0\le r\le s}\big)\Big]=0
\end{equation*}
which completes the proof that $M_t(X)$ is a martingale. Uniqueness of the solution follows from Theorem $4.4.1$ from Ethier and Kurtz's book \cite{EthierKurtz} since the domain of $\SL$ is dense in $\CH^\beta$ for $\beta\in[0,\alpha)$.
\end{proof}

\section{The Brox diffusion on the circle}

Since the potential $\xi=\nabla W$ is periodic with $W(0)=0$, the operator $\SL$ can be considered both on functions on the full line $\IR$ or a finite segment with periodic or Dirichlet boundary condition, we consider the first case here. This can be seen for the operator $\SL_n$ which can be considered on the spaces of smooth periodic functions $C_{\text{per}}^\infty(\IR)$ or smooth compactly supported functions $C_0^\infty(\IR)$. Taking the closure gives two differents operators with respective domains $\CH^2(\IT)$ and $\CH^2(\IR)$. We considered until now the second case, the first one is the content of this section. We consider the operator $\big(\SL,\CD_\Xi(\IT)\big)$
with
\begin{equation*}
\CD_\Xi(\IT)=\{u\in L^2(\IT)\ ;\ u-\PT_{\nabla u}X_1-\PT_{\nabla u}X_2\in\CH^2(\IT)\}
\end{equation*}
and similar form domain. The proof on $\IR$ also gives that $-\SL$ is a $m$-accretive operator in $L^2(\IT,e^{2W}\drm\lambda)$ hence generates a strong Feller semigroup $(e^{-t\SL})_{t\ge0}$ with kernel $p_t(x,y)$ for $t>0$ and $x,y\in\IT$. As for the Laplacian, the spectral properties are very differents in the compact space due to the Sobolev embedding, which is compact. The resolvant $(\SL-c)^{-1}:L^2(\IT)\to\CD_\Xi\subset\CH^\beta$ for any $\beta\in[0,\alpha)$ is compact in $L^2$ and $\SL$ is selfadjoint in $L^2(\IT,e^{2W}\drm\lambda)$ hence it has a pure point spectrum $(\lambda_n)_{n\ge1}$ such that
\begin{equation*}
\lambda_1\ge\lambda_2\ge\ldots\ge\lambda_n\ge\ldots
\end{equation*}
with an associated Hilbert basis of eigenfunctions $(e_n)_{n\ge1}$. The constant function ${\bf 1}=\Gamma{\bf 1}$ belongs to the domain $\CD_\Xi$ with
\begin{equation*}
\SL{\bf 1}=0
\end{equation*}
thus $e_1$ is the constant function and $\lambda_1=0$. The Gaussian lower bounds on the semigroup implies that it is positivity improving in the sense that
\begin{equation*}
f\ge0\quad\implies\quad e^{-t\SL}f>0
\end{equation*}
for $f\in L^2(\IT)\backslash\{0\}$ and any $t>0$. The Krein-Rutman theorem, a generalization of the Perron-Frobenius to infinite dimension, implies that the first eigenvalues is simple, that is
\begin{equation*}
\lambda_1>\lambda_2\ge\lambda_3\ge\ldots
\end{equation*}
and that the associated eigenfunction is of constant sign, for example positive. This agrees with $e_1={\bf 1}$ however one could also consider a finite segment with Dirichlet boundary and the result would still hold. The kernel of $\SL$ is the space of constant functions and this yields an invariant measure on $\IT$ of the process. Since we are working in the Hilbert space $L^2(\IT,e^{2W}\drm\lambda)$, the invariant measure is
\begin{equation*}
\mu(\drm x)=e^{2W(x)}\lambda(\drm x),
\end{equation*}
the measure would be different with Dirichlet boundary conditions. Formaly, we have
\begin{equation*}
\SL^*u=\frac{1}{2}\Delta u-\nabla\big(\nabla W\cdot u\big)=\frac{1}{2}\nabla\big(\nabla u-\nabla 2W\cdot u\big)
\end{equation*}
and $x\mapsto e^{2W(x)}$ is the unique solution to
\begin{equation*}
\nabla f-2\nabla W\cdot f=0
\end{equation*}
on $\IT$ such that $\int_\IT f(x)\drm x=1$. This however only makes sense for the regularized potential $W_n$ and the lower Gaussian bound is crucial to prove the following result as well as the strong Feller property.

\medskip

\begin{theorem}
The measure $\mu$ is the unique invariant measure of the Brox diffusion on the circle.
\end{theorem}

\medskip

\begin{proof}
The measure $\mu$ is indeed invariant for $(e^{-t\SL})_{t\ge0}$ since
\begin{equation*}
\langle\SL u,e^{2W}\rangle=\langle\SL u,{\bf 1}\rangle_{L^2(\IR,e^{2W}\drm\lambda)}=\langle u,\SL{\bf 1}\rangle_{L^2(\IR,e^{2W}\drm\lambda)}
\end{equation*}
using that $\SL$ is symmetric in $L^2(\IT,e^{2W}\drm\lambda)$. Since $\SL{\bf 1}=0$, we get
\begin{equation*}
\langle\SL u,e^{2W}\rangle=0
\end{equation*}
for any $u\in\CD_\Xi$ thus
\begin{equation*}
\langle e^{-t\SL}u,e^{2W}\rangle=\int_\IT\Big(\int_\IT p_t(x,y)u(y)\drm y\Big)e^{2W(x)}\drm x=\int_\IT u(x)e^{2W(x)}\drm x
\end{equation*}
for any $t>0$. Since the domain is dense in $L^2$, we get that $\mu$ is an invariant measure. Because of the Gaussian lower bounds on the semigroup, the associated Markov process is irreductible in the sense that for any $t>0$, we have
\begin{equation*}
\IP(X_{t_0}\in A|X_0=x)>0
\end{equation*}
for any measurable set $A$ with $\lambda(A)>0$ and any $x\in\IT$. This implies that any invariant measure has to be supported on the whole circle $\IT$. Since the set of invariant measures is a convex, the uniqueness of the invariant measure is equivalent to the unicity of extremal invariant measures, that is invariant measures $\nu$ such that for any invariant measures $\nu_1,\nu_1$ with
\begin{equation*}
\eta=\eps\nu_1+(1-\eps)\nu_2
\end{equation*}
and $\eps\in[0,1]$, one has $\eps\in\{0,1\}$. However two extremal invariant measures must be singular since the semigroup is strong Feller, see for example proposition $3.2.5$ from Da Prato and Zabczyk's book \cite{DaPratoZabczyk96}, and the proof is complete.
\end{proof}

A natural question is the convergence of the process $(X_t)_{t\ge0}$ in long time to the measure $\mu$. For any $T>0$, the sequence $(X_{nT})_{n\ge0}$ is a Markov chain with probability transition $e^{-T\SL}$. It is aperiodic and irreductible since of the Gaussian lower bound and the spectral gap implies exponential mixing, that is
\begin{equation*}
\lim_{n\to\infty}\|\SL(X_n)-\mu\|_{\text{TV}}\le C\lambda^n
\end{equation*}
for constants $C>0$ and $0<\lambda<1$ with $\|\cdot\|_{\text{TV}}$ the total variation norm, see for example the very complete book \cite{MeynTweedie09} by Meyn and Tweedie. In particular, the lower bound implies that the circle $\IT$ is a small set in their sense and this is related to the Doeblin criterion. In this continuous setting, this was also investigated and the result still hold, see for example Kontoyiannis and Meyn's work \cite{KontoyiannisMeyn03} and references therein.

\medskip

\begin{corollary}
There exists $C=C(\Xi)>0$ and $\lambda=\lambda(\Xi)\in(0,1)$ such that
\begin{equation*}
\lim_{t\to\infty}\|\SL(X_t)-\mu\|_{\text{TV}}\le C\lambda^t.
\end{equation*}
\end{corollary}

\medskip

It would be natural to consider similar question for the Brox diffusion on the line $\IR$. In the periodic environment, the measure $e^{2W}\drm\lambda$ has infinite mass and can not be normalized to get a probability measure hence the problem should be consider for a two-sided Brownian motion $W$ over $\IR$.

\appendix

\section{Bounds and continuity results}

We give here the needed bounds on $\Phi_n$ and $\Gamma_n$. We also recall the continuity results for the paracontrolled toolkit from \cite{Mouzard}. Note that this is done on a compact manifold but this can be directly adapted to the full space following \cite{BB1}. First, the paraproduct $\P$ and resonant product $\Pi$ allows to decompose a product as
\begin{equation*}
fg=\P_fg+\PI(f,g)+\P_gf,
\end{equation*}
they are built from the heat semigroup $(e^{t\Delta})_{t\ge0}$ as analogue to the Fourier's paraproduct $P$ and resonant product $\Pi$. The paraproduct $\PT$ is intertwined with $\P$ via the relation
\begin{equation*}
\PT\circ\Delta^{-1}=\Delta^{-1}\circ\P
\end{equation*}
and satisfies the same continuity properties as $\P$.

\medskip

\begin{proposition}
Let $\alpha,\beta\in\IR$ be regularity exponents.
\begin{itemize}
	\item[$\centerdot$] If $\alpha\ge0$, then $(f,g)\mapsto\P_fg,\PT_fg$ is continuous from $\CC^\alpha\times\CC^\beta$ to $\CC^\beta$.
	\item[$\centerdot$] If $\alpha<0$, then $(f,g)\mapsto\P_fg,\PT_fg$ is continuous from $\CC^\alpha\times\CC^\beta$ to $\CC^{\alpha+\beta}$.
	\item[$\centerdot$] If $\alpha+\beta>0$, then $(f,g)\mapsto\PI(f,g)$ is continuous from $\CC^\alpha\times\CC^\beta$ to $\CC^{\alpha+\beta}$.
	\item[$\centerdot$] If $\alpha>0$, then $(f,g)\mapsto\P_fg,\PT_fg$ is continuous from $\CC^\alpha\times\CH^\beta$ to $\CH^\beta$ and from $\CH^\alpha\times\CC^\beta$ to $\CH^\beta$.
	\item[$\centerdot$] If $\alpha<0$, then $(f,g)\mapsto\P_fg,\PT_fg$ is continuous from $\CC^\alpha\times\CH^\beta$ to $\CH^{\alpha+\beta}$ and from $\CH^\alpha\times\CC^\beta$ to $\CH^{\alpha+\beta}$.
	\item[$\centerdot$] If $\alpha+\beta>0$, then $(f,g)\mapsto\PI(f,g)$ is continuous from $\CH^\alpha\times\CC^\beta$ to $\CH^{\alpha+\beta}$.
\end{itemize}
\end{proposition}

\medskip

This is enough to prove the bounds on $\Phi_n$ defined by
\begin{equation*}
\Phi_n(u)=u-\PT_{\nabla u}\big(X_1^{(n)}-X_1^{(N_\Xi)}\big)-\PT_{\nabla u}\big(X_2^{(n)}-X_2^{(N_\Xi)}\big)
\end{equation*}
and $\Gamma_n$ by the implict equation
\begin{equation*}
\Gamma_nu^\sharp=\PT_{\nabla\Gamma_nu^\sharp}\big(X_1^{(n)}-X_1^{(N_\Xi)}\big)+\PT_{\nabla\Gamma_nu^\sharp}\big(X_2^{(n)}-X_2^{(N_\Xi)}\big)+u^\sharp.
\end{equation*}
The bounds on $\Phi^{>N}$ and its inverse $\Gamma^{>N}$ are obtained with similar computations and are left to the reader.

\medskip

\begin{proposition}
For any $\beta\in[0,\alpha)$, we have
\begin{equation*}
\|\Phi-\Phi_n\|_{\CH^\beta\to\CH^\beta}\le C_1\|\Xi-\Xi_n\|_{\CX^\alpha}
\end{equation*}
and
\begin{equation*}
\|\textup{Id}-\Gamma\Gamma_n^{-1}\|_{\CH^\beta\to\CH^\beta}\le C_2\|\Xi-\Xi_n\|_{\CX^\alpha}
\end{equation*}
for constants $C_1=C_1(\beta),C_2=C_2(\beta)>0$. In particular, this implies
\begin{equation*}
\|\Gamma-\Gamma_n\|_{\CH^\beta\to\CH^\beta}\le C_3\|\Xi-\Xi_n\|_{\CH^\alpha}
\end{equation*}
for a constants $C_3=C_3(\Xi,\beta)>0$.
\end{proposition}

\medskip

\begin{proof}
The first bound comes directly from
\begin{equation*}
u-\Phi_n(u)=\PT_{\nabla u}\big(X_1^{(n)}-X_1^{(N_\Xi)}\big)+\PT_{\nabla u}\big(X_2^{(n)}-X_2^{(N_\Xi)}\big)
\end{equation*}
and the bounds on $X_i^{(n)}-X_i^{(N_\Xi)}$ for $i\in\{1,2\}$ by the norm of the enhanced noise. For any $u\in\CH^\beta$, we have
\begin{equation*}
u=\Gamma\Gamma^{-1}(u)=\Gamma\Big(u-\PT_{\nabla u}\big(X_1-X_1^{(N_\Xi)}\big)-\PT_{\nabla u}\big(X_2-X_2^{(N_\Xi)}\big)\Big)
\end{equation*}
since $\Gamma^{-1}=\Phi$. We get
\begin{align*}
\|u-\Gamma\Gamma_n^{-1}(u)\|_{\CH^\beta}&=\Big\|\Gamma\Big(u-\PT_{\nabla u}\big(X_1-X_1^{(N_\Xi)}\big)-\PT_{\nabla u}\big(X_2-X_2^{(N_\Xi)}\big)\Big)\\
&\hspace{3cm}-\Gamma\Big(u-\PT_{\nabla u}\big(X_1^{(n)}-X_1^{(N_\Xi)}\big)-\PT_{\nabla u}\big(X_2^{(n)}-X_2^{(N_\Xi)}\big)\Big)\Big\|_{\CH^\beta}\\
&\lesssim\|\PT_{\nabla u}\big(X_1-X_1^{(n)}\big)-\PT_{\nabla u}\big(X_2-X_2^{(n)}\big)\|_{\CH^\beta}.
\end{align*}
As for the bound on $\Phi$, we have for $\beta\in(1,\alpha)$
\begin{equation*}
\|u-\Gamma\Gamma_n^{-1}(u)\|_{\CH^\beta}\lesssim\|\nabla u\|_{L^2}\Big(\|X_1-X_1^{(n)}\|_{\CC^\alpha}+\|X_1-X_2^{(n)}\|_{\CC^{2\alpha-1}}\Big)
\end{equation*}
and for $\beta\in[0,1]$
\begin{equation*}
\|u-\Gamma\Gamma_n^{-1}(u)\|_{\CH^\beta}\lesssim\|\nabla u\|_{\CH^{\beta-1}}\Big(\|X_1-X_1^{(n)}\|_{\CC^\alpha}+\|X_1-X_2^{(n)}\|_{\CC^{2\alpha-1}}\Big)
\end{equation*}
using that $\alpha>1$. This yields
\begin{equation*}
\|u-\Gamma\Gamma_n^{-1}(u)\|_{\CH^\beta}\lesssim\|u\|_{\CH^\beta}\|\Xi-\Xi_n\|_{\CX^\alpha}
\end{equation*}
and the last bound follows with
\begin{equation*}
\|\Gamma-\Gamma_n\|_{\CH^\beta\to\CH^\beta}=\|\big(\textup{Id}-\Gamma\Gamma_n^{-1}\big)\Gamma_n\|_{\CH^\beta\to\CH^\beta}\le\|\textup{Id}-\Gamma\Gamma_n^{-1}\|_{\CH^\beta\to\CH^\beta}\|\Gamma_n\|_{\CH^\beta\to\CH^\beta}.
\end{equation*}
\end{proof}

In this work, we use from \cite{Mouzard} the correctors
\begin{align*}
\DC_\nabla(a_1,a_2,b)&=\PI\big(\nabla\PT_{a_1}a_2,b\big)-a_1\PI(\nabla a_2,b),\\
\DS(a_1,a_2,b)&=\P_b\PT_{a_1}a_2-\P_{a_1}\P_ba_2.
\end{align*}
We state the continuity results in Hölder spaces for simplicity, the needed bound with a mixture of Sobolev and Hölder spaces are also satisfied as for the paraproducts. Note that the proof for $\DC_\nabla$ is similar to one for $\DC$ from \cite{Mouzard}, see \cite{MorinMouzard22} where a similar corrector is used to deal with the first order term.

\medskip

\begin{proposition}
Let $\alpha_1,\alpha_2\in\IR$ and $\beta<0$. Then $(a_1,a_2,b)\mapsto\DS(a_1,a_2,b)$ extends in a unique continuous operator from $\CC^{\alpha_1}\times\CC^{\alpha_2}\times\CC^\beta$ to $\CC^{\alpha_1+\alpha_2+\beta}$.
\end{proposition}

\medskip

\begin{proposition}
Let $\alpha_1\in(0,1)$ and $\alpha_2,\beta\in\IR$. If
\begin{equation*}
\alpha_2+\beta-1<0\quad\text{and}\quad\alpha_1+\alpha_2+\beta-1>0,
\end{equation*}
then $(a_1,a_2,b)\mapsto\DC_\nabla(a_1,a_2,b)$ extends in a unique continuous operator from $\CC^{\alpha_1}\times\CC^{\alpha_2}\times\CC^\beta$ to $\CC^{\alpha_1+\alpha_2+\beta-1}$.
\end{proposition}

\bigskip

\textbf{Acknowledgements :} The author would like to thank Helena Kremp for useful answers about the work \cite{KrempPerkowski22} and Rémi Catellier for an interesting discussion about the diffusion on the circle.

\bibliographystyle{siam}
\bibliography{biblio.bib}

\vspace{2cm}

\noindent \textcolor{gray}{$\bullet$} A. Mouzard --  ENS de Lyon, CNRS, Laboratoire de Physique, F-69342 Lyon, France.\\
{\it E-mail}: antoine.mouzard@math.cnrs.fr

\end{document}